\documentclass[a4paper,12pt]{article}
\usepackage[utf8]{inputenc}
\usepackage[T2A]{fontenc}  
\usepackage[english]{babel}
\usepackage{misccorr}
\usepackage{amssymb}
\usepackage{pgf,tikz}
\usetikzlibrary{arrows}

\title{On periodic approximate solutions of ordinary differential equations}

\author{Wang Shiwei*, A. V. Zorin*, M. A. Konyaeva*, M. D. Malykh*~**, \\
L. A. Sevastianov*~**}
\date{* Patrice Lumumba Peoples' Friendship University of Russia, Moscow, Russia \\
** Joint Institute for Nuclear Research, Dubna, Russia} 
\usepackage{hyperref}

\usepackage[centertags]{amsmath}
\usepackage{amsfonts}
\usepackage{amssymb}
\usepackage{amsthm}
\usepackage{longtable}
\usepackage{multirow}
\usepackage{verbatim}
\usepackage{makeidx}
\makeindex

\def\phi{\varphi}
\def\Phi{\varPhi}
\def\epsilon{\varepsilon}

\theoremstyle{remark}

\theoremstyle{definition}
\newtheorem{Definition}{Definition}

\def\vec{\mathfrak }
\def\dt{\Delta t }

\begin{document}\large \maketitle

\begin{abstract}
The issue of inheriting periodicity of an exact solution of а dynamic system by а difference scheme is considered.
It is shown that some difference schemes (midpoint scheme, Kahan scheme) in some special cases provide
approximate solutions of differential equations, which are periodic sequences. Such solutions are called periodic. A purely algebraic method of finding such solutions is developed. It is shown that midpoint scheme inherits periodicity not only in case of linear oscillator, but also in case of nonlinear oscillator, integrable in elliptic functions.

\end{abstract}

\section{Introduction}

In the study of conservative dynamical systems, periodic solutions play an exceptionally important role. In essence, this is the only class of solutions that is very simply described qualitatively and characterized quantitatively by a single number "--- the period. In those few cases where conservative dynamical systems are explicitly integrated but the solution is not periodic, its qualitative description causes significant difficulties. An example of such a situation is the Kovalevskaya top \cite{Kozlov}. On the contrary, if the problem is not analytically integrable, then particular periodic solutions are used as a kind of ``skeleton'' for describing the entire set of solutions. The most famous example  is the many-body problem. Since Euler, its particular periodic solutions have been sought, and by now numerical methods have been developed that have made it possible to find many thousands of solutions \cite{Marshal, Montgomery-2023}. An absolute breakthrough was the discovery of a figure-eight solution at the end of the last century \cite{Moor, Montgomery-2000,
Montgomery-2001}. The search for periodic solutions to the many-body problem is currently being conducted on supercomputers \cite{Hristov-2022, Hristov-2024}, so more than 5 thousand such solutions were presented at MMCP'2024, found on the HybriLIT platform \cite{Hristov-2024-mmcp}.

In the 1990s, the concept of geometric integrators of differential equations \cite{Lubich} appeared and the question arose of searching for such difference schemes that not only approximate the original dynamic system, but also imitate some properties of its exact solutions. At the end of the last century, the focus was on the inheritance of integrals of motion (conservatism) and Hamiltonian structure (symplecticity) \cite{Lubich}. However, from a practical point of view, a difference scheme can be considered as imitating the original problem if it inherits, in some sense, periodic solutions of the original dynamical problem.

For the case of a linear oscillator, the simplest symplectic scheme, the midpoint scheme, not only preserves the energy integral and the symplectic structure of the problem, but also imitates periodic solutions: the approximate solution in a time step $\dt$ passes  an arc of a circle of fixed length on the phase plane \cite[n.~2]{Malykh-2024-mdpi}. The  difference between the discrete and continuous descriptions is that the lengths of the arcs that the body passes in the discrete and continuous cases differ by a value of the order of $\dt^2$. For nonlinear oscillators integrable in elliptic functions, we were able to trace the inheritance of periodicity for the Kahan scheme \cite[n.~5]{Malykh-2024-mdpi}. The Kahan schemes themselves came into use relatively recently \cite{Sanz-Serna-1994,McLachlan-2013,Suris-2019}.

There are various possible approaches to defining the concept of a periodic approximate solution \cite[n.~5.3]{Malykh-2024-mdpi}. In our opinion, the simplest and at the same time the most general is the idea of a periodic approximate solution as a periodic sequence, proposed in \cite{Malykh-2022-pomi}. The convenience of this definition is that it is naturally formulated for any difference scheme, and the problem of finding periodic solutions becomes purely algebraic.

When starting computer experiments, we assumed that {\it i}) approximate solutions found by the Kahan scheme always inherit the periodicity of the exact solution, and {\it ii}) approximate solutions found by the midpoint scheme and similar ones inherit it only in the linear case. It seemed to us that there should be a direct connection between the reversibility of the scheme, the most important property of the Kahan scheme \cite{Malykh-2024-mdpi}, and the inheritance of periodicity. The experiments presented below refute this hypothesis. 

\section{Approximate periodic solutions}

Let us consider a system of ordinary differential equations
\begin{equation}
\label{eq:ode}
\frac{d\vec x}{dt} = \vec f (\vec x). 
\end{equation}
When using the finite difference method, this system is replaced with a system of algebraic equations
\begin{equation}
\label{eq:schema}
\vec g (\vec x, \hat {\vec x}, \dt)=0,
\end{equation}
relating the values of the variable $\vec x$ at $t$ and at $t+\dt$, the latter is further denoted as $\hat {\vec x}$. From here on, we will consider the step $\dt$ to be a positive number. An approximate solution is a sequence $\vec x_0, \vec x_1, \dots$, whose elements are calculated recursively as solutions of the system
\[
\vec g (\vec x_n, \vec x_{n+1}, \dt)=0
\]
with respect to $\vec x_{n+1}$. 

Here we need to make one important clarification for what follows. If the original problem is linear or if Kahan's scheme is used, then the system \eqref{eq:schema} is linear with respect to $\hat{\vec x}$ and therefore the question of choosing a root when calculating $\vec x_{n+1}$ does not arise. However, in the general case we are faced with the problem of extra roots \cite[n.~3]{Malykh-2024-mdpi}. In numerical analysis this problem is solved as follows. To solve the system with respect to $\hat{\vec x}_{n+1}$, we use an iteration method that converges at sufficiently small steps $\dt$ and returns a root close to $\vec{x}_n$. This rule takes us beyond pure algebra. By analogy with the terminology accepted in elementary mathematics, we will call a recursively determined sequence $\{\vec x_{n}\}$ an approximate solution found according to the scheme \eqref{eq:schema}, regardless of the method of choosing the root. 

In the framework of computer experiments \cite{Malykh-2022-pomi} we noticed that approximate solutions of the Jacobi oscillator found using the Kahan scheme not only approximate periodic solutions, but also for some values of the step $\dt$ form a periodic sequence, i. e., for some $n$ they satisfy the condition
\[
\vec x_{n+m} = \vec x_{m} \quad \forall m=0,1, \dots
\]
In this case, the value $T=n\dt$ does not coincide with the period of the exact solution, but for small $\dt$ it is close to it. In \cite{Malykh-2024-mdpi} we proved this property of the Kahan scheme using its quadrature representation. 

There is nothing to prevent this definition from being extended to the case of arbitrary difference schemes.

\begin{Definition}
A solution found using the scheme \eqref{eq:schema} will be called periodic if it satisfies the condition
\[
\vec x_{n} = \vec x_{0}.
\]
By the period of such a solution we will understand either $n$ or the value $T=n\dt$ depending on the context.
\end{Definition}

Since one point of the solution completely determines the others, the difference scheme \eqref{eq:schema} has a periodic solution if and only if
the system of algebraic equations
\begin{equation}
\label{eq:period}
\vec g (\vec x_0, \vec x_{1}, \dt)=0, \,
\vec g (\vec x_1, \vec x_{2}, \dt)=0
\dots,
\vec g (\vec x_n, \vec x_{0}, \dt)=0
\end{equation}
is compatible. Therefore, the study of periodic approximate solutions always reduces to solving systems of algebraic equations and, at least in theory, is feasible in a finite number of steps. This follows directly from the finiteness of Buchberger's algorithm \cite{CoxLittleOShea}.

It should be noted immediately that for arbitrary values of the step $\dt$ and the initial condition $\vec{x}_0$, the system \eqref{eq:period} containing $n-1$ unknowns $\vec x_1, \dots, \vec x_{n-1}$ and $n+1$ equations is inconsistent. The problem of finding a periodic approximate solution to the Cauchy problem
\begin{equation}
\label{eq:chauchy}
\frac{d\vec x}{dt} = \vec f (\vec x), \quad \vec{x}(0)=\vec{x}_0,
\end{equation}
can reasonably be formulated as follows. For a given $n\in \mathbb{N}$, find a step value $\dt>0$ such that the system \eqref{eq:period} is consistent. 

Generally speaking, the number of unknowns is greater than the number of equations, so periodic solutions should generally exist. However, in practice, only real periods are observed, so the essential question is whether the roots of the \eqref{eq:period} system are real or not.
\section{Computer Experiments}

Using the implementation of Buchberger's algorithm built in the Sage computer algebra system, we investigated the problem of finding a periodic solution to the Cauchy problem for several cases. In the following calculations, instead of the symbolic variable $\dt$, we will use the period $T=n\dt$. In other words, we investigated the consistency of the system
\begin{equation}
\label{eq:eq}
\vec g (\vec x_0, \vec x_{1}, T/n)=0, \,
\vec g (\vec x_1, \vec x_{2}, T/n)=0
\dots,
\vec g (\vec x_n, \vec x_{0}, T/n)=0.
\end{equation}

\subsection{Linear Oscillator}

\begin{figure}
\centering
\includegraphics[scale=0.7]{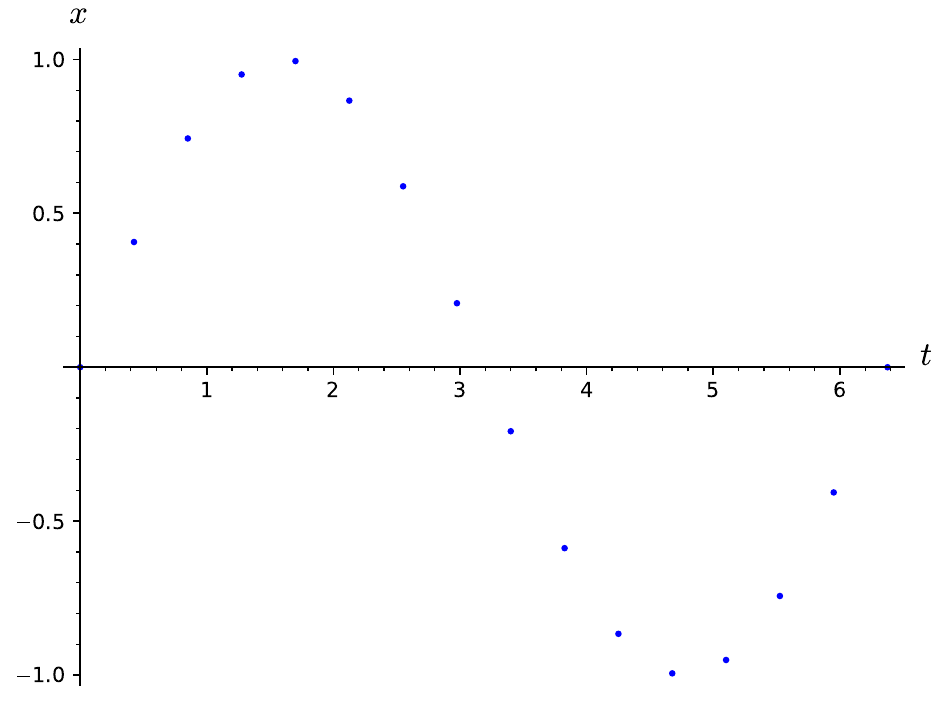}
\caption{Approximate periodic solution of Eq.\eqref{eq:linear} with period $n=15$, found using the midpoint scheme.}
\label{fig:sin}
\end{figure}

First of all, to verify our approach, we considered periodic solutions for the linear oscillator
\begin{equation}
\label{eq:linear}
\frac{dx}{dt}=y, \quad \frac{dy}{dt}=-x, \quad x(0)=0, \quad y(0)=1
\end{equation}
found using the midpoint scheme. From geometric considerations it was previously proved that such solutions exist for all $n\in \mathbb{N}$ and as $n$ tends to infinity the period $T=n\dt \to 2\pi$ \cite{yuying-2019,yuying-2021}. Let us check this.

For \eqref{eq:linear} the system \eqref{eq:eq} is linear with respect to the unknowns $\vec x_1, \dots, \vec x_{n-1}$, their elimination leads to an equation with respect to $T^2$. For the midpoint scheme, this equation has positive roots. For example, for $n=15$ such a root is $T=6.3767 \simeq 2\pi$. At the same time, for $x_0=0, y_0=1$ the plot of the solution (Fig. \ref{fig:sin}) is very similar to the plot of the usual sine. This is quite consistent with the geometric interpretation of the solution of the harmonic oscillator according to the midpoint scheme \cite{yuying-2019,yuying-2021}.

The question of whether the problem \eqref{eq:linear}  has a periodic solution when approximated using the explicit Euler scheme has not been studied previously. Our calculations showed that for $n$ in the region of 10, the equation for $T$ has no real roots different from zero, so the explicit Euler scheme does not yield periodic solutions. From this we can draw the first important conclusion: the existence of exact periodic solutions does not guarantee the existence of approximate periodic solutions.

\begin{figure}
\centering
\includegraphics[scale=0.7]{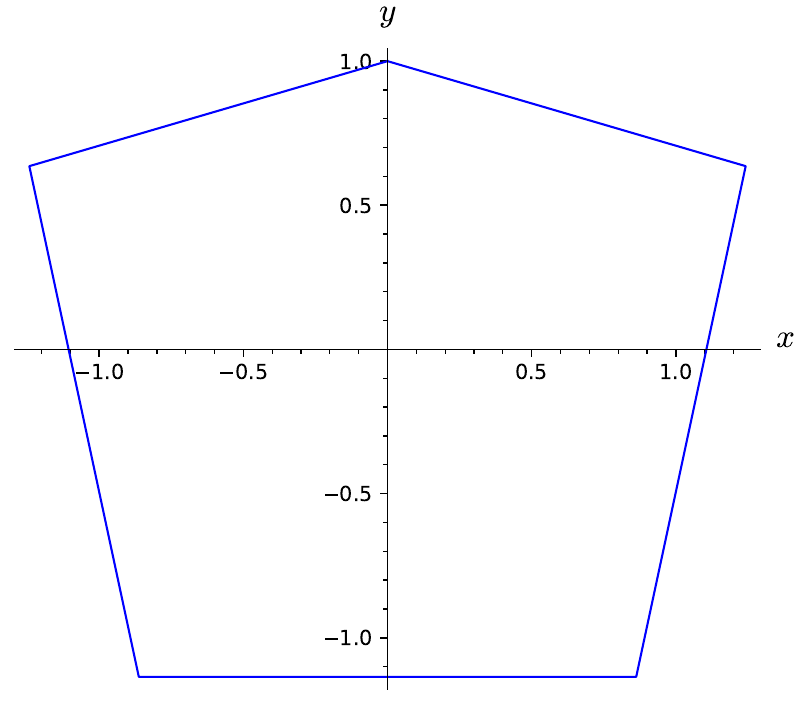}
\includegraphics[scale=0.7]{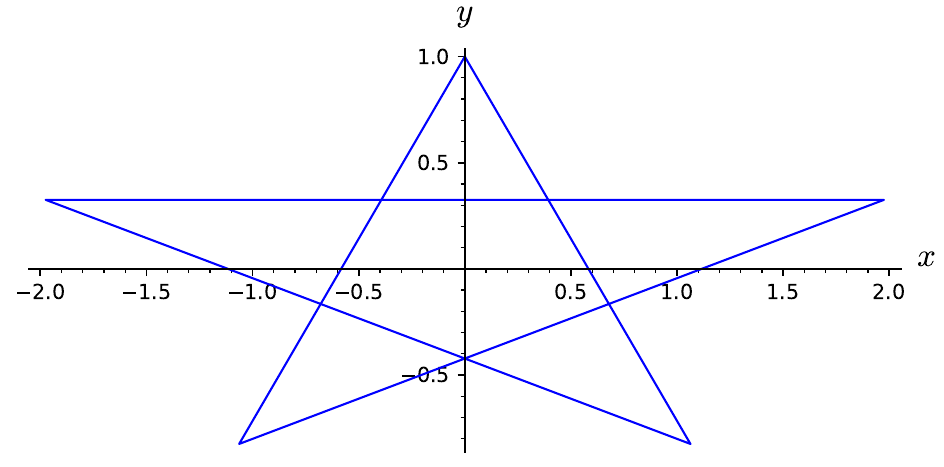}
\caption{Approximate periodic solution of Eq. \eqref{eq:non-linear} with period $n=5$, found using the midpoint scheme at $T=7.59556885597975$ (top) and $T=60.7538695639458$ (bottom).}
\label{fig:pent}
\end{figure}

\subsection{Nonlinear oscillator}

We now turn to the nonlinear oscillator
\begin{equation}
\label{eq:non-linear}
\frac{dx}{dt}=y, \quad \frac{dy}{dt}=-x^3.
\end{equation}
The existence of a periodic solution obtained by the Kahan scheme was pointed out by the authors of Ref. \cite{Malykh-2022-pomi}. Much more interesting is the issue of the existence of periodic solutions found by the midpoint scheme.

For simplicity, let us take a small $n$, namely, $n=5$.
For the midpoint scheme, we obtain a system from which, after eliminating $\vec x_1, \dots, \vec x_4$, we obtain an equation relating the initial point and $T$. This is as it should be, since in the nonlinear case the period always depends on the initial conditions. Let us fix the initial point
\[
\vec x_0=(x_0,y_0)=(0,1),
\]
then we obtain an equation for finding $T^2$, which, in addition to zero, has two positive roots. Both values for $T$ correspond to periodic solutions, shown in Fig.~\ref{fig:pent}. In the first case, a pentagon is obtained in the phase plane, in the second case, a pentagram is obtained.

We observed exactly the same situation for periodic solutions of the elliptic oscillator found using the Kahan scheme \cite{Malykh-2022-pomi}.

It follows that the midpoint scheme imitates the periodicity of the elliptic oscillator \eqref{eq:non-linear} no worse than the Kahan scheme. What is especially interesting here is that the midpoint scheme for nonlinear oscillators specifies a multi-valued correspondence between $\vec x$ and $\hat{\vec x}$. To construct a numerical solution, a root selection rule is added to the \eqref{eq:schema} system; without it, we would get extra roots at each step \cite{Malykh-2024-mdpi}. When compiling the \eqref{eq:eq} system, this rule was by no means taken into account,  nevertheless, this system gave no extra solutions. Thus, it turns out that the periodicity condition $\vec x_{n+m} = \vec x_{m} $ contains the root selection rule.

\subsection{Volterra-Lotka system}


We now turn to the Volterra-Lotka system, taking for definiteness the following numerical values for its parameters:
\begin{equation}
\label{eq:VL}
\left \{ \begin{aligned} &
\frac{d}{dt} x = -x {\left(y - 2\right)} , \quad \frac{d}{dt} y = {\left(x - 2\right)} y , \\ &
x (0)= 1 , \quad y (0)= 2
\end{aligned} \right.
\end{equation}
This problem has a periodic solution, the period being approximately $3.24\dots$. It is easy to find this solution numerically using the Runge-Kutta method with a sufficiently small step.

\begin{figure}
\centering
\includegraphics[scale=0.7]{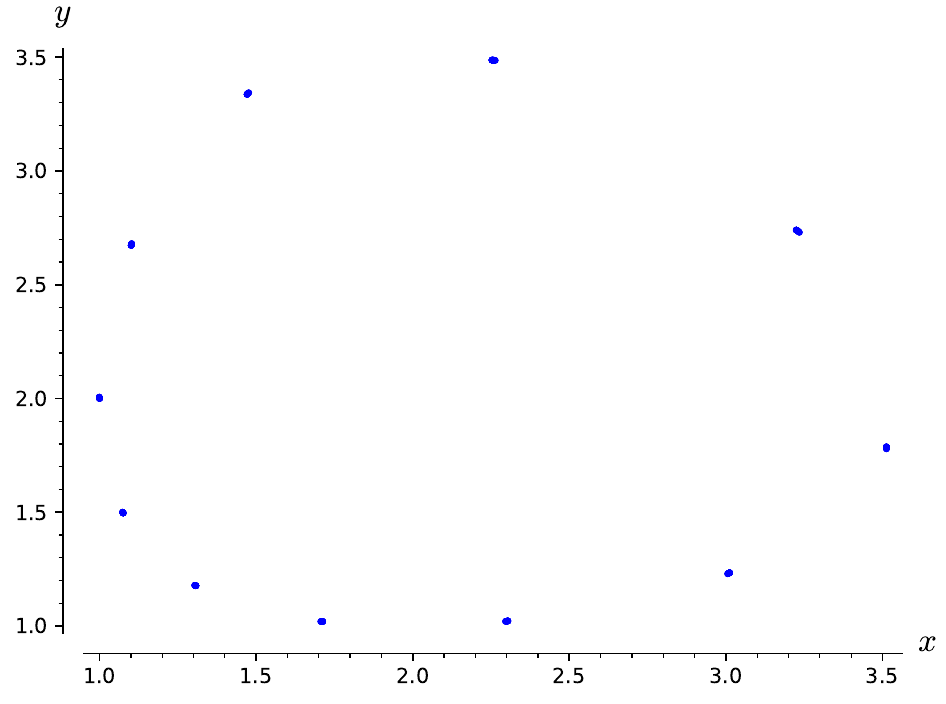}
\caption{Approximate periodic solution of the problem \eqref{eq:VL}, found using the Kahan scheme for $\dt=0.30083\dots$. The first $400$ points are indicated.}
\label{fig:VL}
\end{figure}

Experimenting with Kahan's scheme, we accidentally found that this problem has an approximate solution for $\dt=0.30083\dots$, shown in Fig.~\ref{fig:VL}. It seems that $\vec{x}_{11}=\vec{x}_0$, and small discrepancies in the values are due to rounding error or inaccurate determination of the step $\dt$. Moreover, for the period $T=11\dt$ we get a quite reasonable value of $3.31\dots$.

We investigated the compatibility of system  \eqref{eq:eq} for problem \eqref{eq:VL}. Contrary to our expectations, this system has no nonzero real solutions, including for $n=11$. This means that Kahan's scheme does not always inherit the periodic nature of the solution. Note that in the case of the Volterra-Lotka system, this scheme cannot inherit the integral of motion, which is not algebraic. However, the points of the approximate solution line up along a certain oval. Moreover, for large values of $\dt$, the difference between this oval and the integral curve of Eq. \eqref{eq:VL} is clearly visible.

In parentheses, we note that the study of the compatibility of the system for $n=11$ is a difficult problem for the built-in functions of Sage. We eliminated variables using the fact that the Kahan scheme specifies the Cremona transformations, in the same way that we used in \cite{Malykh-2022-pomi}.

\section{Flipping of an Asymmetric  Top}

In the search for periodic solutions, we get not only an approximate solution, but also its period. This procedure seems to offer us a natural time scale, so that by taking a dozen points, we can, albeit roughly, describe the periodic motion. However, in reality the situation is much more dramatic.

Let us consider Dzhanibekov's top \cite{Petrov-2013} as an example. This system is described by three differential equations
\[
\frac{dp}{dt}=\frac{C-B}{A} qr, \quad
\frac{dq}{dt}=\frac{A-C}{B}pr, \quad
\frac{dr}{dt}=\frac{B-A}{C}pq,
\]
where $A,B,C$ are the principal moments of inertia, and is integrable in elliptic functions \cite{Golubev-1953}. The case described by Dzhanibekov is interesting due to the sharp flips of the top, which present an obvious difficulty in numerical integration. In this case, the proportion of the principal moments of inertia is
$A:B:C=8:7:2$. Then, if we spin the top about the middle (unstable) axis of inertia and slightly disturb the initial conditions, that is, take
$p(0)=0$, $q(0)=1$ and $r(0)=\epsilon$, where $\epsilon$ is a small parameter, then the top will be flipping: half the period $q \simeq 1$, the second half of the period $q\simeq -1$. The transition between these values occurs the faster, the smaller $\epsilon$, and the period itself grows as $\ln \epsilon^{-1}$ \cite{Petrov-2013}.

\begin{figure}
\centering
\includegraphics[scale=0.7]{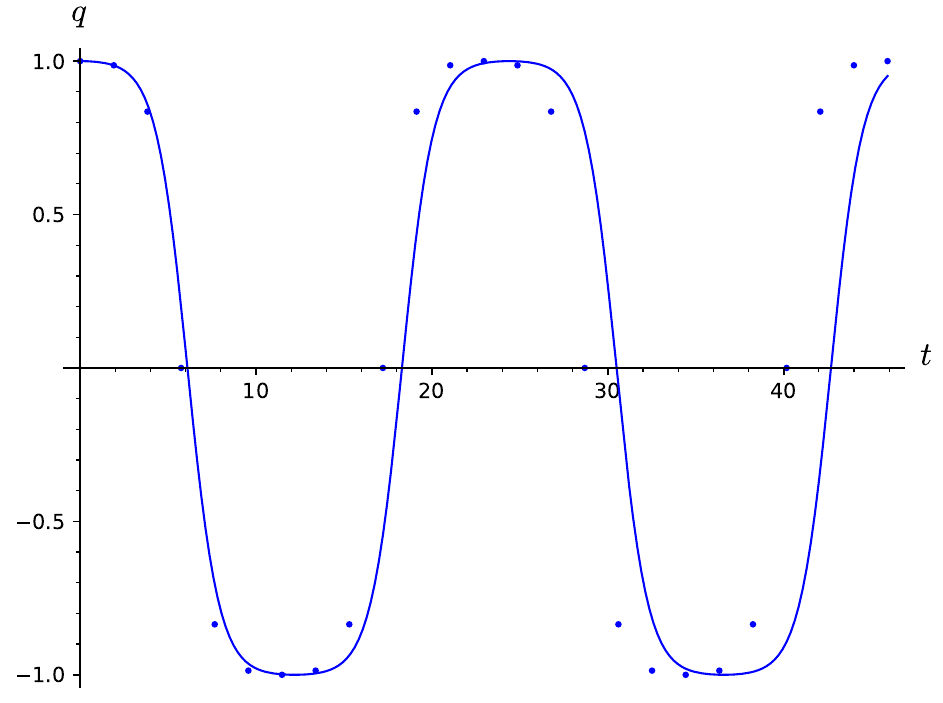}
\includegraphics[scale=0.7]{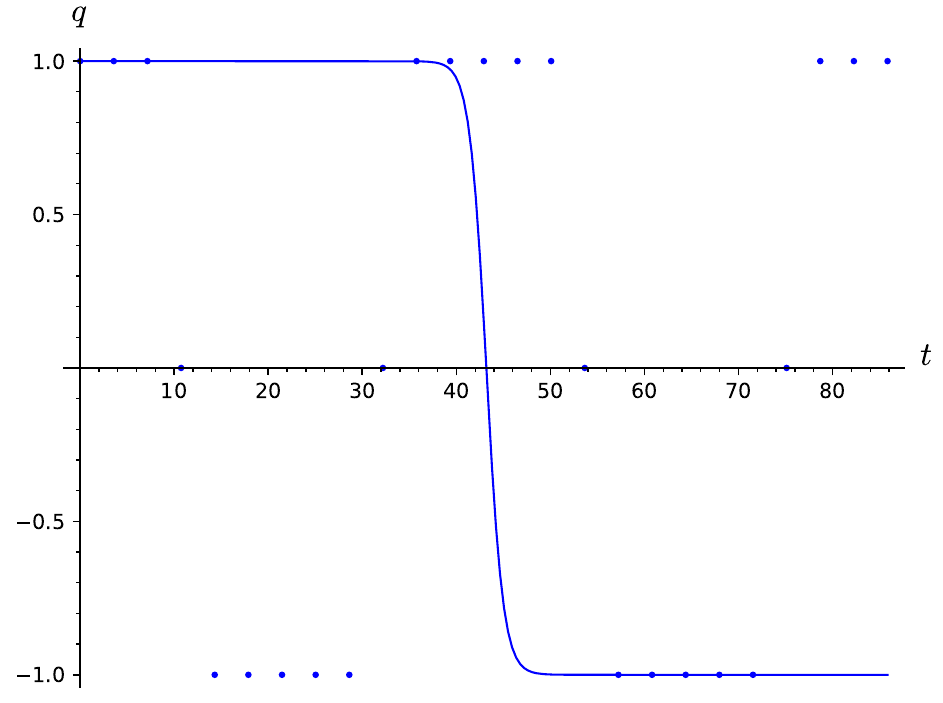}
\caption{Exact and approximate periodic solutions describing the top flips for $\epsilon = 10^{-1}$ (top) and $\epsilon = 10^{-10}$ (bottom).}
\label{fig:ABC}
\end{figure}

The relationship between the exact solution in elliptic functions and the approximate solution according to Kahan's scheme for the top is described in Ref. \cite[n.~5]{Malykh-2024-mdpi}. Let us now calculate the periodic solution for $n=12$. To make calculations faster, as in the many-body problem, we consider half the period rather than the entire period,  taking $p=r=0, q=1$ as the initial point, and $p=r=0, q=-1$ as the final point. For $\epsilon=10^{-1}$, we obtained three periodic solutions with periods
$22.952$, $42.569$ and $42.932$. The smallest of them corresponds to a periodic solution that perfectly imitates the exact one, see Fig. \ref{fig:ABC}. Only 12 points per period in the second-order scheme ensure good accuracy in determining the period (the error less than 5\%) and perfect qualitative repetition of the motion.

However, with a decrease in the parameter $\epsilon$, the situation drastically changes  quantitatively although still unchanged qualitatively. For example, for $\epsilon=10^{-10}$, the solution with the smallest period of $42.925$ is qualitatively similar to the exact one, but the value of this period differs by several times from the period of the exact solution, see Fig. \ref{fig:ABC}. Obviously, the transition from $q=1$ to $q=-1$ is described by one intermediate point $q=0$. This means that the accuracy of determining the period depends on the number of points per transition layer ratgheer than on the number of points per period.

Thus, with flipping top, we encounter a very unusual situation. Qualitatively, the approximate solution coincides with the exact one: the top really flips periodically, the transition layer is sharply outlined. However, the error in determining the period is quite noticeable.

\section{Conclusion}

Let us summarize the experiments. 
First of all, we found  that approximate solutions determined using the Kahan scheme do not always inherit the periodicity of the exact solution. The inheritance occurs for linear and elliptic oscillators and does not occur for example, for the Volterra-Lotka system.

Secondly, approximate solutions found by the midpoint scheme inherit the periodicity of the exact solution not only in the case of linear oscillators, but also in the case of elliptic oscillators. At the same time, difference schemes that introduce dissipation or anti-dissipation into the model, for example, the explicit Euler scheme, cannot yield periodic solutions. Therefore, we can say that the class of difference schemes that give periodic solutions for a given dynamic system is wider than the class of Kahan schemes. The question of how much wider it is requires additional study.

We have left aside the undoubtedly important question of the numerical solution of the system of equations \eqref{eq:eq}. In the case where several dozen points must fit into a period, the analytical solution of the system becomes very expensive.
Our experiments with the Volterra-Lotka system showed that a solution to this system with a small residual can exist even when the system is actually incompatible.

Within the framework of the proposed approach to finding periodic solutions, we specify the number of solution points per period, and obtain not only a periodic solution, but also its period. Unfortunately, the estimate obtained in this way for the period of the exact solution may differ from the period of the exact solution quite significantly. The key obstacle to the closeness of these periods is a rapid change in the exact solution over short periods of time, for example, sharp flips of a top.

The proposed approach to finding periodic solutions is equally applicable to explicit and implicit difference schemes. The periodicity condition in a sense ``regularizes'' an implicit scheme, removing unnecessary roots. This circumstance is equally unexpected and curious. The point is that unnecessary roots are the main problem of all implicit schemes, due to which calculations by implicit schemes become many times more complicated than calculations by explicit schemes. The results of our experiments allow us to hope that this problem can be circumvented when calculating periodic solutions.

\bibliographystyle{ugost2008}
\bibliography{kahan}

\end{document}